\documentclass{amsart}
\usepackage{amssymb,latexsym}
\usepackage{graphics}
\usepackage[cp1251]{inputenc}
\usepackage[english,russian]{babel}
\usepackage{amscd}
\usepackage{amsmath}
\usepackage[matrix,arrow,curve]{xy}
\usepackage{euscript}
\usepackage{lscape}
\numberwithin{equation}{section}

\theoremstyle{plain}
\newtheorem{prop}{Proposition}

\theoremstyle{definition}
\newtheorem{lemma}{Lemma}

\newtheorem{definition}{Definition}
\theoremstyle{remark}
\newtheorem{remark}{Remark}

\numberwithin{corollary}{section}
\newcommand{\Os}{\mathcal{O}}
\newcommand{\Gt}{\mathbb{G}}
\newcommand{\Zt}{\mathbb{Z}}
\newcommand{\Pt}{\mathbb{P}}
\newcommand{\Pcat}{\mathcal{P}}
\newcommand{\rar}{\rightarrow}

\newcommand{\Gm}{\mathbb{G}_m}

\begin{document}
\centerline{\bf{ALGEBRAIC ANALOGUE OF ATIYAH'S THEOREM.}}
\bigskip
\centerline{\Small{ALISA KNIZEL AND ALEXANDER NESHITOV.}}
\bigskip
{\centerline{\parbox{11cm}{\footnotesize{{\bf{Abstract}} In topology there is a well-known theorem of Atyiah which states that
for a connected compact Lie group G there is an isomorphism $\widehat{R(G)}\cong K_0(BG)$  where BG is the classifying space of G.
In the present paper we consider an algebraic analogue of this theorem.
In the paper by B.Totaro [8]
it is shown that the  $\varprojlim K_0(BG_i)$ is equal $\widehat{R(G)}$
for a specially chosen sequence $BG_i.$
However, to compute $K_0(BG)$ one needs to prove that $\varprojlim^1 K_1(BG_i)$ vanishes.
For a split reductive group $G$ over a field we present another aproach and prove 
that the Borel construction induces a ring isomorphism
$\widehat{R(G)}_{I_G}\cong K_0(BG),$ where $I_G$ is the fundamental ideal of $R(G).$ This statement can be 
formulated in terms of equivariant $K$-theory as:
$\widehat{K_0^G(pt)}_{I_G}\cong K_0(BG).$
The main aim of the present paper is to extend this result for higher K-groups. Namely, we prove that there is a natural isomorphism
$$\widehat{K_n^G(pt)}_{I_G}\cong K_n(BG),$$
where $K_n^G(pt)$ is Thomason's $G$-equivariant K-theory defined in[3], $BG$ is a motivic \'etale classifying space introduced by Voevodsky and Morel in [6] and for a motivic space $X$ the group $K_n(X)$ is defined as $Hom_{H_{\bullet}}(S^n_s\wedge(X_+),BGL\times\mathbb{Z})$([10], thm. 6.5)
}}}}
\section{Introduction}
{\parbox{13cm}{

We will work over a field $k.$ Morel and Voevodsky in ([6],Def. 4.2.4, Prop 4.2.6) constructed the \'etale classifying space of a linear algebraic group $G$ in the form $BG = \bigcup BG_{m},$ where $BG_{m}= EG_{m}/G$ and $EG_{m}$ are $k$-smooth algebraic varieties with a free $G$-action, connected by a sequence of G-equivariant closed embeddings $i_{k}$ }} 
 $$ \begin{CD}
 \cdots @>{i_{m-1}} >> (EG)_{m} @>{i_{m}} >> (EG)_{m+1} @>{i_{m+1}} >> \cdots 
 \end{CD}
 $$

{\parbox{13cm}{ The motivic space $EG=\bigcup EG_{m}$ is $\mathbb{A}^{1}$ -- contractible with a free $G$-action([6], Prop. 4.2.3). 
We consider a split reductive affine algebraic group $G.$ A $G$-equivariant vector bundle over the $pt=Spec(k)$ is the same
as a $k$-rational representation of $G.$ So, we will identify these two categories. Notice that this identification respects
the tensor products. In particular, we will identify Thomason's $K_0^G(pt)$ with the representation ring of $k$-rational representations 
$R(G)$ of the group $G.$

The Borel construction sends a $G$-equivariant vector bundle $V$ over the point to the vector bundle $V_{m}=(V\times EG_{m})/G$ over $BG_m.$
This construction respects tensor products. Therefore it induces a $K_0^G(pt)$-modules map
$\phi_m:K_n^G(pt)\rar K_n(BG_m).$ Obviously, $\phi_m=\bar{i_m}^*\circ \phi_{m+1},$ where $\bar{i_m}:BG_m\rar BG_{m+1}$ is an embedding
induced by $i_m.$ As we will prove below (Proposition 3), $K_n(BG)=\varprojlim K_{n}(BG_m).$ Combining all these, we get an $K_0^G(pt)$-module map
$$ \Phi_n: K_n^G(pt)\rar K_n(BG).$$ 
We will write $Borel^G_n$ for $\Phi_n.$ Let $I_G$ be the kernel of the augmentation $K_0^G(pt)\rar K_0(pt)=\mathbb{Z}.$\\

{\bf{Theorem 1 (main)}} :\\
In the following diagram both maps are $K_0^G(pt)$-module isomorphisms
$$\xymatrix{\widehat{K_n^G(pt)_{I_G}}\ar[rr]^{\widehat{Borel^G_n}} && \widehat{K_n(BG)}_{I_G} && K_n(BG)\ar[ll]_{completion_G} \\},$$
where $\widehat{Borel_n^G}$ is the $I_G$ completion of $Borel^G_n,$ and $completion_G$ is the canonical map.}

{\parbox{13cm}{The main idea of the proof is the reduction to a Borel subgroup B of G.
For the Borel subgroup $B$ the $K_0^B(pt)$-modules $K_n(BB)$ and $K_n^B(pt)$ can be computed explicitly. It results in}}
\begin{flushleft}
{{\bf{Theorem 2}}{ The Borel construction induces an isomorphism}}
\end{flushleft}
$$\xymatrix{\widehat{K_n^B(pt)}_{I_B}\ar[rr]^{\widehat{Borel^B_n}} && \widehat{K_n(BB)}_{I_B} && K_n(BB)\ar[ll]_{\cong} \\}$$
To make a reduction to the Theorem 2 we prove
\begin{flushleft}
{\bf{Theorem 3}}{ There is a commutative diagram of the form:}
\end{flushleft}

\begin{equation}
\xymatrix{
     \widehat{K_n^G(pt)}_{I_G}\ar[rr]^{\widehat{Borel^G_n}}\ar[d]_{res} && \widehat{K_n(BG)}_{I_G}\ar[d]_{\widehat{p^*}} && K_n(BG)\ar[ll]_{}\ar[d]_{p^*}\\
     \widehat{K_n^B(pt)}_{I_B}\ar[rr]^{\widehat{Borel^B_n}}\ar[d]_{ind} && \widehat{K_n(BB)}_{I_B}\ar[d]_{\widehat{p_*}} && K_n(BB)\ar[ll]_{\cong}\ar[d]_{p_*}\\
     \widehat{K_n^G(pt)}_{I_G}\ar[rr]^{\widehat{Borel^G_n}} && \widehat{K_n(BG)}_{I_G} && K_n(BG) \ar[ll]_{} \\
}
\label{main reduction diagram}
\end{equation}
With $ind\circ res=id,$ $\widehat{p_*}\circ\widehat{p^*}=id,$ $p_*\circ p^*=id.$\\

{\bf Remark.} Clearly main theorem follows from theorem 2 and theorem 3.
Now we are working on the generalization of this result for the case of non-split reductive algebraic groups.\\ 

{\bf{Acknowledgements.}} Authors are grateful to Prof. Ivan Panin for constant attention and useful suggestions concerning the subject of this paper.
\section{Auxiliary results}

{\begin{flushleft}{In this section we prove some properties of pullback and pushforward morphisms for~$K_n^G$~functor.
Thomason in [3] developed G-equivariant K-theory.(c.f. Merkurjev's paper [2])}\end{flushleft}}
\definition{Let X be a G-variety. We consider an action $\mu_x:G\times X \rightarrow X$ and a projection
$p_x:G\times X\rightarrow X.$ Let M be an $\Os_X$-module. Following [2] we will call $M$ a $G$-module 
if there is an isomorphism of $\Os_{G\times X}$-modules $\alpha:\mu_X^*(M)\rightarrow p_X^*(M)$ such that the cocycle condition holds:
$$p_{23}^*(\alpha)\circ(id_G\times \mu_x)^*(\alpha)=(m\times id_X)^*(\alpha)$$
where $p_{23}:G\times G\times X\rightarrow G\times X$ is a projection and $m:G\times G\rightarrow G$ is a product morphism.
}
\smallskip
\definition{We denote by $\Pcat(G;X)$ a category of locally free $G$-modules on $X.$}
\definition{Equivariant $K$-functor $K_n^G(X)$ define as $K_n(\Pcat(G;X))$}

\lemma{Let $f:X\rightarrow Y$ be an equivariant morphism and let $M$ be a G-module on Y

Then $f^*M$ has a structure of G-module on X.\\
Proof:\\
Consider the following diagram:


$$\xymatrix{
     G\times X \ar@<1ex>[d]^{\mu_X}\ar@<-1ex>[d]_{p_X} \ar[rr]^{id_G\times f} && G\times Y\ar@<1ex>[d]^{\mu_Y}\ar@<-1ex>[d]_{p_Y}  \\
     X\ar[rr]^{f}        && Y        \\
}
$$
We construct $\alpha$ as a composition of isomorphisms:
$$\xymatrix{
     p_X^*f^*M && (id_G\times f)^*p_Y^*M\ar[ll]^{\cong}  \\
     \mu_X^*f^*M\ar[rr]^{\cong}\ar[u]^{\alpha}        && (id_G\times f)^*\mu_Y^*\ar[u]^{(id_G\times f)^*\beta}        \\
}
$$
Here $\beta$ is a $G$-module structure on $M$. The cocylce condition for $\alpha$ immediately follows from the cocycle condition for $\beta.$\\
{\bf{Corollary 1.}} For any equivariant $f:X\rightarrow Y$ we have exact functor 
$$f^*:\Pcat(G;Y)\rightarrow \Pcat(G;X).$$
It induces a pullback morphism $K_n^G(Y)\rightarrow K_n^G(X).$ To simplify notation we will also denote it by $f^*.$
\remark
Let X,Y be G-varieties, f be a G-morphism, M be an $\Os_X$-module, N be an $\Os_{G\times X}$-module,  F=$id_G\times f.$ Consider the diagram:
$$\xymatrix{
     G\times G\times X \ar@<-1ex>[d]_{m\times id_X}\ar@<1ex>[d]^{id_G\times \mu_X} \ar[r]^{id_G\times F} & G\times G\times Y\ar@<-1ex>[d]_{m\times id_Y}\ar@<1ex>[d]^{id_G\times \mu_Y} && G\times G\times X \ar@<-1ex>[d]_{m\times id_X}\ar@<1ex>[d]^{p_{23_X}} \ar[r]^{id_G\times F} & G\times G\times Y\ar@<-1ex>[d]_{m\times id_Y}\ar@<1ex>[d]^{p_{23_Y}} \\
     G\times X \ar[d]^{\mu_X} \ar[r]^{F} & G\times Y\ar[d]^{\mu_Y} && G\times X \ar[d]^{p_X} \ar[r]^{F} & G\times Y\ar[d]^{p_Y} \\
     X\ar[r]^{f}        & Y && X\ar[r]^{f}        & Y        \\
}
$$

{\bf{Notation 1:}}\\
Since all vertical arrows are flat, we have natural isomorphisms ([1] Prop. 9.3):
$hh_{\mu}(M):\mu_Y^*R^if_*M\rightarrow R^iF_*\mu_X^*M;$ \\
$hh_{p}(M):p_Y^*R^if_*M\rightarrow R^iF_*p_X^*M $\\
$hh_{m\times id}(N):(m\times id_Y)^*R^iF_*N\rightarrow R^i(id_G\times F)_*(m\times id_X)^*N;$\\ $hh_{p_{23}}(N):p_{23_Y}^*R^iF_*N\rightarrow R^i(id_G\times F)_*p_{23_X}^*N; $\\
$hh_{id_G\times \mu}(N):(id_G\times \mu_Y)^*R^iF_*N\rightarrow R^i(id_G\times F)_*(id_G\times \mu_X)^*N;$\\

Note that since $\mu_Y\circ(m \times id_Y)=\mu_Y\circ(id_G\times\mu_Y),$ two  isomorphisms coincide:\\
$hh_{\mu,id_G\times\mu}(M):(id\times\mu_Y)^*\mu_Y^*R^if_*M\rightarrow R^i(id_G\times F)_*(id_G\times\mu_X)^*\mu_X^*M$ and \\
$hh_{\mu,m\times  id}(M):(m\times id_Y)^*\mu_Y^*R^if_*M\rightarrow R^i(id_G\times F)_*(m\times id_X)^*\mu_X^*M$\\
Similarly, there is another pair of equal isomorphisms:\\
$hh_{p,p_{23}}(M):p_{23_Y}^*p_Y^*R^if_*M\rightarrow R^i(id_G\times F)_*p_{23_X}^*p_X^*M$ and \\
$hh_{p,m\times  id}(M):(m\times id_Y)^*p_Y^*R^if_*M\rightarrow R^i(id_G\times F)_*(m\times id_X)^*p_X^*M$\\
\medskip
We need the following lemma about compostion of this isomorphisms.
\lemma{ Consider the following diagram:
$$\xymatrix{
     X_3\ar[rr]^{f_3}\ar[d]_{T}  && Y_3\ar[d]_{Q}\\
     X_2\ar[rr]^{f_2}\ar[d]_{t}  && Y_2\ar[d]_{q}\\
     X_1\ar[rr]^{f_1}  && Y_1\\
}
$$
Here $q$ and $Q$ are flat, $X_2=X_1\times_{Y_1}Y_2,$ $X_3=X_2\times_{Y_2}Y_3.$ Let M be an $\Os_{X_1}$-module.
Define \\
$hh_1:q^*R^i f_{1*}\rightarrow R^i f_{2*}t^* $\\
$hh_{12}:Q^*q^*R^i f_{1*}\rightarrow R^i f_{3*}T^*t^*$\\
$hh_2:Q^*R^i f_{2*}\rightarrow R^i f_{3*}T^*$ to be natural isomorphisms given by Prop. 9.3 [1].
Then the following diagram commutes:
$$\xymatrix{
    Q^*q^*R^if_{1*}M \ar[rd]_{hh_{12}(M)} \ar[rr]^{Q^*hh_1(M)} &&  Q^*R^if_{2*}t^*M\ar[ld]_{hh_2(t^*M)}\\
                        & R^if_{3*}T^*t^*M
}
$$
Proof:\\
Since the statement is local on $Y_i,$ we consider the case when all $Y_i$ are affine, $Y_i=Spec A_i.$
If $F$ is $R$- module, we will denote by $\widetilde{F}$ the corresponding sheaf on $Spec \ R.$ 
Recall the construction of $hh_1.$
Let $M$ be an $\Os_{X_1}$-module. Then\\
 $R^if_*(M)=\widetilde{H^i(X_1,M)};$
$q^*R^if_{1*}M=\widetilde{A_2\otimes_{A_1}H^i(X_1,M)};$
$R^if_{2*}t^*M=\widetilde{H^i(X_2,t^*M)}.$\\
Let $U_i$ be an affine covering of $X_1.$ Denote by $K=\check{C}(X_1,M)$ the corresponding \v{C}hech complex.
Since $Y_1$ and $Y_2$ are affine, $t^{-1}(U_i)$ is the affine covering of $X_2.$ For this covering we have that 
$A_2\otimes_{A_1}K$ is a \v{C}hech complex of $X_2$-module $t^*M.$ Then $hh_1$ is an obvious morphism 
$$A_2\otimes_{A_1}H^i(K)\rightarrow H^i(A_2\otimes_{A_1}K)$$ which becomes an isomorphism since $A_2$ is flat over $A_1.$
In similar way one can construct $hh_{12}$ and $hh_2.$ Then one can rewrite the diagram as 
$$\xymatrix{
    {A_3\otimes_{A_2}A_2\otimes_{A_1}H^i(K)}\ar[rd]_{hh_{12}(M)} \ar[rr]^{id\otimes hh_1} && {A_3\otimes_{A_2}H^i(A_2\otimes_{A_1}K)}\ar[ld]_{hh_2(t^*M)} \\
                        & {H^i(A_3\otimes_{A_1}K)}
}
$$
Which is trivially commutative.
}

\lemma{Let $f:X\rightarrow Y$ be an equivariant morphism and $M$ be a G-module on X.
Then for any $i$  $R^if_*M$ has a structure of G-module on Y.\\
Proof:\\
Let $\beta:\mu_x^* M\longrightarrow p_X^* M$ be the $G$-structure on M.
Consider the following base-change diagram:

$$\xymatrix{
     G\times X \ar@<0.5ex>[rr]^{\mu_X}\ar@<-0.5ex>[rr]_{p_X} \ar[d]^{id_G\times f} && X\ar[d]^{f}  \\
     G\times Y \ar@<0.5ex>[rr]^{\mu_Y}\ar@<-0.5ex>[rr]_{p_Y}        && Y        \\
}
$$
Since $\mu_Y$ and $p_Y$ are flat, we use Proposition 9.3 from [1]. Sheaf isomorphisms $hh_{\mu}(M)$ and $hh_p(M)$ are described in Notation 1.
Define $\alpha$ to be the unique isomorphism such that the following diagram commutes :

 $$\xymatrix{
     \mu_Y^*R^if_*M\ar[rr]^{hh_{\mu}(M)}\ar[d]^{\alpha} && R^i(id\times f)_*\mu_X^*M\ar[d]^{R^i(id\times f)_*\beta}  \\
     p_Y^*R^if_*M\ar[rr]^{hh_p(M)}        && R^i(id\times f)_*p_X^*M\\
}
$$
Now we have to check the cocycle condition for $\alpha:$ \\
$p_{23}*(\alpha)\circ(id_G\times \mu_Y)^*(\alpha)=(m\times id_Y)^*(\alpha)$\\
This means commutativity of this diagfam:
 $$\xymatrix{
     (id_G\times\mu_Y)^*\mu_Y^*R^if_*M\ar[rrr]^{p_{23}*(\alpha)\circ(id_G\times \mu_Y)^*(\alpha)}\ar@{=}[d] &&& p_{23}^*p_Y^*R^if_*M\ar@{=}[d] \\
     (m\times id_Y)^*\mu_Y^*R^if_*M\ar[rrr]^{(m\times id_Y)^*(\alpha)}        &&& (m\times id_Y)^*p_Y^*R^if_*M\\
}
$$
Let $F=id_G\times f.$ Subdivide this diagram into the following blocks:

 $$\xymatrix{
     (id_G\times\mu_Y)^*\mu_Y^*R^if_*M\ar[rr]_{\displaystyle 1}\ar[d]^{\cong}_{hh_{\mu,id\times\mu}} && p_{23_Y}^*p_Y^*R^if_*M\ar[d]^{\cong}_{hh_{p,p_{23}}}  \\
     R^i(id_G\times F)_*(id_G\times \mu_X)^*\mu_X^*M\ar[rr]_{\displaystyle 2}\ar@{=}[d]        && R^i(id_G\times F)_*p_{23_X}^*p_X^*M\ar@{=}[d]\\
     R^i(id_G\times F)_*(m\times id_X)^*\mu_X^*M\ar[rr]_{\displaystyle 3}\ar[d]^{\cong}_{hh_{m\times id}^{-1}(\mu_X^*M)} && R^i(id_G\times F)(m\times id_X)^*p_X^*M\ar[d]^{\cong}_{hh_{m\times id}^{-1}(p_X^*M)}\\
     (m\times id_Y)^*R^i(id_G\times f)_*\mu_X^*M\ar[rr]_{\displaystyle 4}\ar[d]^{\cong}_{(m\times id_Y)^*hh_{\mu}^{-1}(M)}        && (m\times id_Y)^*R^i(id_G\times f)_*p_X^*M\ar[d]^{\cong}_{(m\times id_Y)^*hh_{p}^{-1}(M)}\\
     (m\times id_Y)^*\mu_Y^*R^if_*M\ar[rr]^{(m\times id_Y)^*(\alpha)}        && (m\times id_Y)^*p_Y^*R^if_*M\\
     }
$$
Square 2 is an image of cocycle diagram for M and therefore commutative.\\
Square 3 arises from functor isomorphism \\
$R^i(id_G\times id_G \times f)_*(m\times id_X)^*\cong (m\times id_Y)^*R^i(id_G\times f)_*$ ([1], Prop. 9.3) applied to 
G-module structure $\beta:\mu_X^* M\rightarrow p_X^* M$ So, it commutes.\\
Square 4 is commutative by definition of $\alpha.$\\
It remains to show the commutativity of square 1. Let $\widetilde{F}=id_G\times F.$  Rewrite square 1 as follows:
$$\footnotesize{\xymatrix@M=0pt{
     (id\times\mu_Y)^*\mu^*R^if_*M\ar[r]_{\displaystyle 1.1}\ar[d]^{\cong}_{(id\times \mu)^*hh_{\mu}} &(id\times\mu_Y)^*p_Y^*R^if_*M \ar@{=}[r]_{\displaystyle 1.2}\ar[d]^{\cong}_{(id_G\times \mu_Y)^*hh_{p}} &p_{23_Y}^*\mu_Y^*R^if_*M\ar[r]_{\displaystyle 1.3}\ar[d]^{\cong}_{p_{23_Y}^*hh_{\mu}} &p_{23_Y}^*p_Y^*R^if_*M\ar[d]^{\cong}_{p_{23_Y}^*hh_{p}}  \\
    (id\times\mu_Y)^*R^iF_*\mu_X^*M\ar[r]_{\displaystyle 1.4}\ar[d]^{\cong}_{hh_{id\times\mu}(\mu^*M)} &(id\times \mu_Y)R^iF_*p_X^*M\ar[d]^{\cong}_{hh_{id\times\mu}(p_X^*M)} &p_{23_Y}^*R^iF_*\mu_X^*M\ar[r]_{\displaystyle 1.5}\ar[d]^{\cong}_{hh_{p_{23}}(\mu_X^*M)} &p_{23_Y}^*R^iF_*p_X^*M\ar[d]^{\cong}_{hh_{p_{23}}(p_X^*M)}\\
    R^i\widetilde{F}_*(id\times \mu_X)^*\mu_X^*M\ar[r] &R^i\widetilde{F}_*(id\times \mu_X)^*p_X^*M\ar@{=}[r] &R^i\widetilde{F}_*p_{23_X}^*\mu_X^*M\ar[r] &R^i\widetilde{F}_*p_{23_X}p_X^*M\\
}}
$$
Square 1.1 is an image of functor $(id_G\times \mu_Y)^*$ applied to the diagram that defines $\alpha.$ Thus it is commutative.
Commutativity of 1.2 follows from Lemma 2 and Prop. 9.3 [1] applied to the base-change diagram
$$\xymatrix{
     G\times G\times X \ar[rr]^{id_G\times id_G\times f}\ar[d]_{p_Y\circ(id\times \mu_Y)=\mu_Y\circ p_{23_Y}} && G\times G\times Y\ar[d]^{p_Y\circ(id\times \mu_X)=\mu_X\circ p_{23_X}}  \\
     X \ar[rr]^{f}      && Y        \\
}
$$
Square 1.3 is an image of functor $p_{23_Y}^*$ applied to the diagram defining $\alpha$ and therefore commutes.\\
Prop. 9.3 [1] gives us an isomorphsm of functors $(id_G\times \mu_Y)^*R^iF_*\cong R^i\widetilde{F}(id_G\times \mu_X)^*.$
Applying this isomorphism to $\beta:\mu_X^*M\rightarrow p_X^*M$ we get commutativity of the square 1.4.\\
In a similar way we get commutativity of the square 1.5.\\
So commutativity of 1-4 is proved. According to Lemma 2 the composition of vertical arrows is the identity. So $\alpha$ satisfies the cocycle condition.

}

{\bf{Corollary 2.}}{  If f is projective we can define the pushforward map $f_*:K_0^G(X)\rightarrow K_0^G(Y)$ by sending M to  the alternating sum of $R^if_*(M).$}\\
To describe the pushforward for higher $K$-functors we need the following lemmas:
\lemma{(Equivariant version of Proposition 9.3[1].) Consider the base change diagram 
$$\xymatrix{
     A\ar[rr]^{F}\ar[d]_{Q}  && B\ar[d]_{q}\\
     X\ar[rr]^{f}  && Y\\
}
$$
where $X,Y,A,B$ are $G$-varieties; $f,F,Q,q$ are $G$-morphisms; $f$ is flat.\\
Let M be a $G$-module on $B.$ Then there is  a natural $G$-module isomorphism on $X:$
$$f^*R^iq_*M\rightarrow R^iQ_*F^*M.$$
Proof:\\
By Propostion 9.3 from [1] we have a natural isomorphsim of $\Os_X$-modules\\ $hh_{X,Y,A,B}:f^*R^iq_*M\rightarrow R^iQ_*F^*M.$ We need to check that $hh_{X,Y,A,B}$ is a $G$-morphism. That means commutativity of the following diagram:
$$\xymatrix{
     \mu_X^*f^*R^iq_*M\ar[d]^{\mu_X^*hh_{X,Y,A,B}}\ar[rr]_{G-structure} &&  p_X^*f^*R^iq_*M\ar[d]^{p_X^*hh_{X,Y,A,B}} \\
     \mu_X^*R^iQ_*F^*M\ar[rr]_{G-structure} && p_X^*R^iQ_*F^*M   \\
}
$$
Consider the diagram:
$$\xymatrix{
G\times A \ar[rr]^{id\times F}\ar[dd]_{id\times Q}\ar@<0.5ex>^{p_A}[dr]\ar@<-0.5ex>_{\mu_A}[dr]
& & G\times B \ar'[d]_{id\times q}[dd]\ar@<0.5ex>^{p_B}[dr]\ar@<-0.5ex>_{\mu_B}[dr]
\\
& A \ar[rr]_<<<<<<<<<<<{F}\ar[dd]_<<<<<<<{Q}
& & B \ar[dd]_{q}
\\
G\times X \ar'[r][rr]_{id\times f}\ar@<0.5ex>^{p_X}[dr]\ar@<-0.5ex>_{\mu_X}[dr]
& & G\times Y \ar@<0.5ex>^{p_Y}[dr]\ar@<-0.5ex>_{\mu_Y}[dr]
\\
& X \ar[rr]^{f}
& & Y 
}
$$
For any square in this cube denote by $hh$ (with corresponding subscript) the isomorphism arising from prop. 9.3[1], applied to this square.
We rewrite the $G$-structure diagram:\\
 
 $$\xymatrix{
     \mu_X^*f^*R^iq_*M\ar@{=}[d]\ar[rr]_{\displaystyle 1} &&  p_X^*f^*R^iq_*M\ar@{=}[d] \\
     (id\times f)^*\mu_Y^*R^iq_*M\ar[rr]_{\displaystyle 2}\ar[d]_{(id\times f)^*hh_{G\times Y,Y,G\times B,B}^{\mu}(M)} && (id\times f)^*p_Y^*R^iq_*M\ar[d]_{(id\times f)^*hh_{G\times Y,Y,G\times B,B}^{p}(M)}   \\
     (id\times f)^*R^i(id\times q)_*\mu_B^*M\ar[rr]_{\displaystyle 3}\ar[d]_{hh_{G\times X, G\times Y, G\times A, G\times B}(\mu_B^*M)} && (id\times f)^*R^i(id\times q)_*p_B^*M\ar[d]_{hh_{G\times X, G\times Y, G\times A, G\times B}(p_B^*M)}\\
     R^i(id\times Q)_*(id\times F)^*\mu_B^*M\ar@{=}[d]\ar[rr]_{\displaystyle 4} && R^i(id\times Q)_*(id\times F)^*p_B^*M\ar@{=}[d]\\
R^i(id\times Q)_*\mu_A^*F^*M\ar[d]_{hh_{G\times X,X,G\times A,A}^{\mu}}\ar[rr]_{\displaystyle 5 } && R^i(id\times Q)_*p_A^*F^*M\ar[d]_{hh_{G\times X,X,G\times A,A}^{p}}\\
     \mu_X^*R^iQ_*F^*M\ar[rr] && p_X^*R^iQ_*F^*M\\
}
$$
Square 1 is commutative because of the definition of the $G$-structure on pullback.\\
Square 2 is an $(id\times f)^*$ image of the $G$-structure diagram for $R^iq_*M.$ Thus it commutes.\\
Square 3 arises from the functor isomorphism $(id\times f)^*R^i(id\times q)_*\rightarrow R^i(id\times Q)_*(id\times F)^*$ applied to
the $G$-structure isomorphism $\mu_B^*M\rightarrow p_B^*M.$ So it commutes.\\
Square 4 is commutative because of the definition of the $G$-structure on pullback.\\
Square 5 is commutative by the definition of the $G$-structure on $R^iQ_*F^*M.$\\
By lemma 2 compositions of vertical arrows are equal to $\mu_X^*hh_{X,Y,A,B}$ and $p_X^*hh_{X,Y,A,B}.$
This concludes the proof of Lemma 4.
}

\lemma{Let $X,Y$ be smooth G-varieties, G - a smooth reductive affine algebraic group and $\pi:X\times Y\rightarrow Y$ a projection. Moreover let $X$ be projective and $Y$ be connected \\
Denote by $\Pcat_{\pi}(G; X\times Y)$ the full subcategory of $\Pcat(G; X\times Y)$ consisting of locally free $G$-modules $P$ such that $R^k\pi_*P=0$ for $k>0.$\\
Then any $G$- module $M$ possesses a finite length resolution of the form\\
$M\rightarrow P^0\rightarrow P^1\rightarrow\ldots\rightarrow P^N\rightarrow 0$ with $P^i \in OB(\Pcat_{\pi}(G;X\times Y))$\\
Proof:\\
First we prove that for every $M$ there is an embedding $M\hookrightarrow P^0.$ We will construct $P^0$ in the form of $M(n)$ for a large enough $n.$
To do this, we construct a very ample $G$-equivariant sheaf $\Os_X(1)$ and an 
$G$-equivariant embedding $i:X\hookrightarrow \mathbb{P}^n$ such that $\Os_X(1)=i^*\Os_{\mathbb{P}}(1).$
Let $L$ be a very ample line bundle. By corollary 1.6 of [5]  $L^{\otimes k}$ is
G-equivariant for some $k$. Then it defines the action of $G$ on $V=\Gamma(X,L^{\otimes k})$ and equivariant morphism $i:X\rightarrow \mathbb{P}(V)$ which is an embedding since $L^{\otimes k}$ is very ample. Then we set $\Os_X(1)=L^{\otimes k}.$\\
The standard embedding of the tautological bundle $\tau_{\mathbb{P}(V)}\hookrightarrow V\times\mathbb{P}(V)$ gives us a $G$-equivariant embedding of locally free sheaves $\Os_{\mathbb{P}(V)}(-1)\hookrightarrow \Os_{\mathbb{P}(V)}\oplus\ldots\oplus \Os_{\mathbb{P}(V)}.$ After twisting by $\Os_{\mathbb{P}}(1)$ we have 
$\Os_{\mathbb{P}(V)}\hookrightarrow \Os_{\mathbb{P}(V)}(1)\oplus\ldots\oplus \Os_{\mathbb{P}(V)}(1).$ Inductively we have the $G$-equivariant embedding $\Os_{\mathbb{P}(V)}\hookrightarrow \Os_{\mathbb{P}(V)}(n)\oplus\ldots\oplus \Os_{\mathbb{P}(V)}(n).$
Applying $i^*$ we get
$$\Os_X\hookrightarrow \Os_X(n)\oplus\ldots\oplus \Os_X(n).$$
Define $\Os_{X\times Y}(1)=\pi^*\Os_{X}(1).$ Applying $\pi^*$ we get an equivariant embedding
$$M\hookrightarrow M(n)\oplus\ldots\oplus M(n).$$
for an arbitrary locally free $G$-module M. Clearely it's cokernel is $G$-equivariant. It's easy to check that it is a locally free sheaf.
Then for every locally free $G$-module there is a resolution consisisting of direct sums of modules of the form $M(n)$
\\
Let us show that $M(n)$ lies in $\Pcat_{\pi}(G;X\times Y)$ for a large enough $n.$
$R^k\pi_*M(n)$ is associated to a presheaf $V\mapsto H^k(X\times V,M(n)).$ Consider a finite affine covering $V_i$ of Y. 
By Serre's theorem $H^k(X\times V_i,M(n))$ equals zero for $n>n_i$. Thus $R^k\pi_*M(n)=0$ for $n>n_M=max\{n_i\}.$ 

It remains
to show that this resolution ends at some finite step.
Let $N=dim X\times Y.$
Let $C^0$ be a cokernel or the first resolution step:
$0\rightarrow M\rightarrow P^0\rightarrow C^0\rightarrow 0.$ Then we have the exact sequence
$$0=R^N\pi_*P^0\rightarrow R^N\pi_*C^0\rightarrow R^{N+1}\pi_*M=0.$$ So, $R^N\pi_*C^0=0.$
For the second cokernel $C^1$
we have the exact sequence $0\rar C^0\rar P^1 \rar C^1 \rar 0.$ Then
$$0=R^{N-1}\pi_*P^1\rightarrow R^{N-1}\pi_*C^1\rightarrow R^{N}\pi_*C^0=0.$$
So, $R^{N-1}\pi_*C^{N-1}=0.$
By induction we have all $R^k\pi_*C^N=0.$ Then $C^N\in Ob(\Pcat_{\pi}(G;X\times Y)).$ 
}\\

{\bf{Corollary 3.}}{  Let f be a $G$-equivariant projective $f:X\rightarrow Y.$ That means, there is a $G$-equivariant decomposition
$$
\xymatrix{
& Y\times \Pt^n\ar[d]^{\pi_Y}\\
X \ar@{^{(}->}[ur]^{i}\ar[r]^{f} & Y \\
}
$$
Here $i$ is a closed embedding and $\pi_Y$ a projection. Since all $R^k i_*M=0$ for any $G$-module $M,$ $k>0,$
we have two exact functors $i_*\Pcat(G;X)\rightarrow \Pcat(G;Y\times\Pt^n)$
and $\pi_{Y*}:\Pcat_{\pi_Y}(G;Y\times \Pt^n)\rightarrow \Pcat(G;Y)$
By Quillen's theorem, the inclusion $\Pcat_{\pi_Y}(G;Y\times \Pt^n)\subseteq \Pcat(G;Y\times \Pt^n)$
induces an isomorphism $$\xymatrix{
K_n(\Pcat_{\pi_Y}(G;Y\times \Pt^n))\ar[rr]^{\alpha} && K_n(\Pcat(G;Y\times \Pt^n))=K_n^G(Y\times\Pt^n)\\
}$$
Then we can describe the pushforward map 
$f_*:K_n^G(X)\rightarrow K_n^G(Y)$ as the following compostion:
$$\xymatrix{ 
K_n^G(X)\ar[r]^{K_n(i_*)} & K_n^G(Y\times \Pt^n) \ar[r]^{\alpha^{-1}} & K_n(\Pcat_{\pi_Y}(G;Y\times \Pt^n)) \ar[r] & K_n(\Pcat(G;Y))=K_n^G(Y) \\ 
}$$
\smallskip
\lemma{Under the notation of Lemma 5, we have a commutative up to an isomorphism diagram of exact functors.}

\begin{equation}
\xymatrix{
\Pcat_{\pi_{EG_j}}(G;EG_j\times G/B)\ar[d]_{\pi_{EG_j *}} && \Pcat_{\pi_{EG_{j+1}}}(G;EG_{j+1}\times G/B)\ar[ll]^{(i_j\times id)^*}\ar[d]_{\pi_{EG_{j+1 *}}}\\
\Pcat(G;EG_j) && \Pcat(G;EG_{j+1})\ar[ll]^{i_j^*}\\
}\label{Categories diagram pushforward}
\end{equation}
Proof:\\
To simplfy notation let $\pi_j=\pi_{EG_j}$ and $\Pcat_j=\Pcat_{\pi_{EG_j}}(G;EG_j\times G/B)$
Let us prove that $\mathcal{P}_{j+1}$ is mapped to $\mathcal{P}_{j}$ under $(i_{j} \times id)^*$
Let $M\in Ob(\mathcal{P}_{j+1}).$ Let $dim(EG_{j}\times G/B)=N.$ Then $R^{N+1}\pi_{j *}(i_j\times id)^*M=0.$ By corollary 2\S 5 of [4]}

$$R^N\pi_{j *}(i_j\times id)^*M\otimes_{\Os_{EG_j}}k(y)=H^N(EG_j\times \{y\},(i_j\times id)^*M)=H^N(EG_j\times \{y\},M)=0 $$\\
Then $ R^N\pi_{j *}(i_j\times id)^*M=0.$ By induction we obtain that all $R^k\pi_{j *}i_j^*M=0$ for $k>0.$ Then $i_j^*M\in Ob(\Pcat).$
Now we prove the commutativity of the diagram \ref{Categories diagram pushforward} up to a natural isomorphism.
By remark 9.3.1 of [1] we have a natural morphism $hh:i_j^*\pi_{j+1 *}M\rightarrow \pi_{j *} (i_j\times id)^* M.$
One can easily see that for any $y \in EG_j$  the following diagram commutes:
$$\xymatrix{
\pi_{j *}(i_j\times id)^*M\otimes k(y)\ar[d]_{(1)} && i_j^*\pi_{j+1 *}M \otimes k(y)\ar[ll]^{hh\otimes k(y)}\ar@{=}[d] \\
\Gamma(y\times G/B,(i_j\times id)^*M)\ar@{=}[d] && \pi_{j+1 *}M\otimes k(y)\ar[d]_{(2)}\\
\Gamma(y\times G/B,M) && \Gamma(y\times G/B,M)\ar@{=}[ll]\\
}.$$
Here the arrows (1) and (2) are natural isomorphisms given by corollary 2 \S 5 of [4].
So, $hh\otimes k(y)$ is an isomorphism for any point $y$ of $EG_j.$ Therefore $hh$ is a natural isomorphism. So, the diagram (\ref{Categories diagram pushforward}) is commutative.

\lemma{Under the notation of Lemma 5, for each $j\geqslant 0$ the functor $$\pi_j^*:\Pcat(G;EG_j)\rightarrow \Pcat(G;EG_j\times G/B)$$ takes values in the subcategory $\Pcat_{\pi_{EG_j}}(G;EG_j\times G/B).$ As a consequence the following diagram of exact functors commutes up to a natural isomorphism.
\begin{equation}
\xymatrix{
\Pcat(G;EG_j)\ar[d]_{\pi_j^*} && \Pcat(G;EG_{j+1})\ar[ll]^{i_j^*}\ar[d]_{\pi_{j+1}^*}\\
\Pcat_{\pi_{EG_j}}(G;EG_j\times G/B) && \Pcat_{\pi_{EG_{j+1}}}(G;EG_{j+1}\times G/B)\ar[ll]^{(i_j\times id)^*}\\
}\label{Categories diagram pullback}
\end{equation}
Proof:\\
To simplfy notation let $\pi_j=\pi_{EG_j}$ and $\Pcat_j=\Pcat_{\pi_{EG_j}}(G;EG_j\times G/B)$
First we prove that $\pi_j^*$ maps $\Pcat(G;EG_j)$ to $\Pcat_j.$
Let $M$ be an object of $\Pcat(G;EG_j).$ Then $R^k\pi_{j *}\pi_j^*M$ is associated to the presheaf $V\mapsto H^k(V\times G/B,\pi_j^*M).$
Let $V$ be an affine open subset of $EG_j$.
Let $\{U_n\}$ be an affine covering of $G/B.$ For any intersection $W=U_{n_1}\cap\ldots\cap U_{n_k}.$ we have
$$\pi_j^*M(V\times W)=M(V)\otimes_{\Os_{EG_j}(V)}\Os_{EG_j\times G/B}(V\times W)=M(V)\otimes_k \Os_{G/B}(W).$$
Then \v{C}hech complex $\check{C}(\{V\times U_n\},\pi_j^*M)$ equals $M(V)\otimes_k\check{C}(\{U_n\},\Os_{G/B}).$ 
Consequently, $H^k(V\times G/B,\pi_j^*M)=M(V)\otimes_k H^k(G/B,\Os_{G/B}).$\\
By proposition 4.5 from [7] $H^k(G/B,\Os_{G/B})=0$ for $k>0.$ Then $\pi_{j *}M \in Ob (\Pcat_j).$
The commutativity of (\ref{Categories diagram pullback}) trivally follows from the equality $\pi_{j+1}\circ (i_j\times id)=i_j\circ\pi_j.$

\lemma{Composition $\pi_{EG_j *}\circ\pi_{EG_j}^*$ is naturally isomorphic to $id_{\Pcat(G;EG_j)}$:}
$$\xymatrix{
\Pcat(G;EG_j)\ar[r]^<<<<<<{\pi_{EG_j}^*} & \Pcat_{\pi_{EG_j}}(EG_j\times G/B)\ar[r]^<<<<<<{\pi_{EG_j *}} & \Pcat(G;EG_j)\\
}$$
Proof:\\
Let $M\in Ob(\Pcat(G;EG_j))$ The sheaf
$\pi_{EG_j *}\pi_{EG_j}^*M$ is associated to presheaf $V\mapsto \pi_{EG_j}^*(M)(V\times G/B).$
Since $\pi_{EG_j}^*M$ is a sheaf associated to $W\mapsto M(\pi_{EG_j}(W))$ we see that $\pi_{EG_j *}\pi_{EG_j}^*M$ is associated to the presheaf $V\mapsto M(V).$
So, in category of presheaves $\pi_{EG_j *}\pi_{EG_j}^*\cong id.$ Applying the sheaffication functor to this isomorphism, we get a natural isomorphism $\pi_{EG_j *}\pi_{EG_j}^*M\cong M.$\\
The same reasoning proofs 

\lemma{Composition $\pi_{pt *}\circ\pi_{pt}^*$ is naturally isomorphic to $id_{\Pcat(G;pt)}$:}
$$\xymatrix{
\Pcat(G;pt)\ar[r]^<<<<<<{\pi_{pt}^*} & \Pcat_{\pi_{pt}}(G/B)\ar[r]^<<<<<<{\pi_{pt *}} & \Pcat(G;pt)\\
}$$

\lemma{Under the notation of Lemma 5, we have a commutative up to an isomorphism diagram of exact functors.}

\begin{equation}
\xymatrix{
\Pcat_{\pi_{pt}}(G;G/B)\ar[d]_{\pi_{pt *}}\ar[rr]^{\pi_{G/B}^*} && \Pcat_{\pi_{EG_{j}}}(G;EG_{j}\times G/B)\ar[d]_{\pi_{EG_{j *}}}\\
\Pcat(G;pt)\ar[rr]^{\pi_{pt}^*} && \Pcat(G;EG_{j})\\
}\label{Categories diagram pushforward2}
\end{equation}
Proof:\\
Let us prove that $\Pcat_{\pi_{pt}}(G;G/B)$ is mapped to $\Pcat_{\pi_{EG_{j}}}(G;EG_{j}\times G/B)$ under $\pi_{G/B}^*.$
Let $M\in Ob(\Pcat_{\pi_{pt}}(G;G/B)).$ Then by prop 9.3[1] $R^k\pi_{EG_j}(\pi_{G/B}^*M)$ is isomorphic to $\pi_{pt}^*(R^k\pi_{pt*}M).$
The latter sheaf is zero by definition of $\Pcat_{\pi_{pt}}(G;G/B)$ for $k>0.$
So, for $k>0$ we have $R^k\pi_{EG_j}(\pi_{G/B}^*M)=0$ then $\pi_{G/B}^*M \in Ob(\Pcat_{\pi_{EG_{j}}}(G;EG_{j}\times G/B)).$

Commutativity of the diagram \ref{Categories diagram pushforward2} follows immediately from lemma 4.

\lemma{Under the notation of Lemma 5, functor $$\pi_{pt}^*:\Pcat(G;pt)\rightarrow \Pcat(G;G/B)$$ takes values in the subcategory $\Pcat_{\pi_{pt}}(G; G/B).$ As a consequence the following diagram of exact functors commutes up to a natural isomorphism.
\begin{equation}
\xymatrix{
\Pcat(G;pt)\ar[d]_{\pi_{pt}^*}\ar[rr]^{\pi_{pt}^*} && \Pcat(G;EG_{j})\ar[d]_{\pi_{EG_j}^*}\\
\Pcat_{\pi_{pt}}(G;G/B)\ar[rr]^{\pi_{G/B}^*} && \Pcat_{\pi_{EG_{j}}}(G;EG_{j}\times G/B)\\
}\label{Categories diagram pullback2}
\end{equation}
Proof:\\
We prove that $\pi_{pt}^*$ maps $\Pcat(G;pt)$ to $\Pcat_{\pi_{pt}^*}(G;G/B).$
Let $M$ be an object of $\Pcat(G;pt).$ Then $R^k\pi_{pt *}\pi_{pt}^*M$ is a vector space $H^k(G/B,\pi_{pt}^*M).$
Let $\{U_n\}$ be an affine covering of $G/B.$ For any intersection $W=U_{n_1}\cap\ldots\cap U_{n_k}.$ we have
$$\pi_{pt}^*M(W)=M\otimes_{k}\Os_{G/B}(W).$$
Then \v{C}hech complex $\check{C}(\{U_n\},\pi_{pt}^*M)$ equals $M\otimes_k\check{C}(\{U_n\},\Os_{G/B}).$ 
Consequently, $H^k(G/B,\pi_j^*M)=M\otimes_k H^k(G/B,\Os_{G/B}).$\\
By proposition 4.5 from [7] $H^k(G/B,\Os_{G/B})=0$ for $k>0.$ Then $\pi_{pt *}M \in Ob(\Pcat_{\pi_{pt}}(G;G/B)).$
The commutativity of (\ref{Categories diagram pullback2}) trivally follows from the equality $\pi_{pt}\circ\pi_{EG_j}=\pi_{pt}\circ\pi_{G/B}.$\\

{\bf{Remark 1}} As we can see from proofs of lemmas 6-11, we can replace $G/B$ by any projective $G$-variety $X$ such that $h^0(X,\Os_X)=1$ and $h^i(X,\Os_X)=0.$ for $i>0.$

\bigskip
{\prop{ There is a commutative diagram with $\pi_{EG_i *}\circ \pi_{EG_i}^*=id_{K_n^G(EG_i)},$ $\pi_{pt *}\pi_{pt}^*=id_{K_n^G(pt)}$}}
\\
$$\xymatrix{
     K_n^G(pt) \ar[d]_{\pi_{pt}^*}\ar[rr]^{\pi_{pt}^*}  && K_n^G(EG_i)\ar[d]_{\pi_{EG_i}^*}            \\
     K_n^G(G/B)\ar[rr]^{\pi_{G/B}^*}\ar[d]_{\pi_{pt}*}  && K_n^G(EG_i\times G/B)\ar[d]_{\pi_{EG_i *}}  \\
     K_n^G{pt}\ar[rr] ^{\pi_{pt}^*}        && K_n^G(EG_i)                  \\
}
$$
Proof:\\
By lemmas 10 and 11 we get the following commutative categories diagram with exact arrows
\begin{equation}
\xymatrix{
\Pcat(G;pt)\ar[d]_{\pi_{pt}^*}\ar[rr]^{\pi_{pt}^*} && \Pcat(G;EG_{j})\ar[d]_{\pi_{EG_j}^*}\\
\Pcat_{\pi_{pt}}(G;G/B)\ar[d]_{\pi_{pt *}}\ar[rr]^{\pi_{G/B}^*} && \Pcat_{\pi_{EG_j}}(EG_j\times G/B)\ar[d]_{\pi_{EG_j *}}\\
\Pcat(G;pt)\ar[rr]^{\pi_{pt}^*} && \Pcat(G;EG_{j})\\
}\label{finn}
\end{equation}
Recall that, by Quillen's theorem and lemma 5, categories inclusion $\Pcat_{\pi_{pt}}(G;G/B)\subseteq \Pcat(G;G/B)$
induces an isomorphism $K_n(\Pcat(G;G/B))\rightarrow K_n(\Pcat_{\pi_{pt}}(G;G/B)).$
Then, applying $K_n$ to diagram (\ref{finn}) gives us
$$\xymatrix{
     K_n^G(pt) \ar[d]_{\pi_{pt}^*}\ar[rr]^{\pi_{pt}^*}  && K_n^G(EG_i)\ar[d]_{\pi_{EG_i}^*}            \\
     K_n^G(G/B)\ar[rr]^{\pi_{G/B}^*}\ar[d]_{\pi_{pt}*}  && K_n^G(EG_i\times G/B)\ar[d]_{\pi_{EG_i *}}  \\
     K_n^G{pt}\ar[rr] ^{\pi_{pt}^*}        && K_n^G(EG_i)                  \\
}
$$
Equalities $\pi_{EG_i *}\circ \pi_{EG_i}^*=id_{K_n^G(EG_i)}$ and $\pi_{pt *}\pi_{pt}^*=id_{K_n^G(pt)}$ immediately follow from lemma 8 and 9. 

{\bf{Remark 2}} In particular, we get a well-known fact that the natural ring map $R(G)\rightarrow R(B)$ is injective.\\

{\bf{Remark 3}} By remark 1, we can replace $G/B$ in statement of proposition 1 by any projective $G$-variety $X$ such that $h^0(X,\Os_X)=1$ and $h^i(X,\Os_X)=0.$ for $i>0.$\\
{\prop{The $I_B$-adic topology of $R(B)$ coincides with the $I_G\cdot R(B)$-adic topology.}}\\
Proof:\\
Let $T$ be a maximal torus in $G.$ Then $R(B)=R(T)$ and $I_B=I_T,$ where $I_T$ is the ideal of zero-dimensional representations of $T.$ 
We will prove that $\sqrt{I_G\cdot R(T)}=I_T.$ Denote by $W=N_G(T)/T$ the Weil group of $G.$
The group $W$ acts by conjugation on $R(T).$ It is known that $W$ is a finite group and $R(G)$ is the ring of invariants of $W$:
$R(G)=R(T)^{W}.$ 
We prove the following statement:\\ \smallskip
If $q$ is a prime ideal of $R(T)$ and $q\cap R(G)\supseteq I_G.$ Then $q\supseteq I_T.$\\
\smallskip
Let $x\in I_T.$ Let $n=|W|$ and $W=\{\sigma_1,\ldots,\sigma_n\}.$
For any symmetric polynomial $f$ we have that $f(x^{\sigma_1}\ldots x^{\sigma_n})$ is invariant under $W$-action.
Then $f(x^{\sigma_1}\ldots x^{\sigma_n})\in R(G)\cap I_T=I_G\subseteq R(G)\cap q.$ Then $f(x^{\sigma_1}\ldots x^{\sigma_n})\in q.$
Denote by $f_1\ldots f_n$ the elementary symmetric polinomials. It is easy to see that $x$ is a root of polynomial
$$\prod\limits_{i=1}^{n}(t-x^{\sigma_i})=t^n-f_1(x^{\sigma_1}\ldots x^{\sigma_n})t^{n-1}+\ldots +(-1)^nf_n(x^{\sigma_1}\ldots x^{\sigma_n}).$$
So we have $x^n-f_1(x^{\sigma_1}\ldots x^{\sigma_n})x^{n-1}+\ldots +(-1)^nf_n(x^{\sigma_1}\ldots x^{\sigma_n})=0.$ \\
Then $x^n=-(-f_1(x^{\sigma_1}\ldots x^{\sigma_n})x^{n-1}+\ldots +(-1)^nf_n(x^{\sigma_1}\ldots x^{\sigma_n}))\in q.$
So $x^n\in q.$ Since $q$ is prime, $x\in q.$ This ends the proof of the statement.\\
Consider $A=\{p\mid p-prime, p\supseteq I_G\cdot R(T)\}$
Our statement implies that $I_T$ is a minimal element of $A.$  So,$$\sqrt{I_G\cdot R(T)}=\bigcap\limits_{p\in A}=I_T.$$
Since $R(B)=R(T)$ and $I_B=I_T,$ we get $\sqrt{I_G\cdot R(B)}=I_B.$ Since $R(B)$ is noetherian, it implies that $I_B^m\subseteq I_G\cdot R(B)$ for some m. Then $I_B$ and $I_G\cdot R(B)$ determine the same topology on $R(B).$
{\prop{$K_n(BG)=\varprojlim K_n(BG_i)$}}\\
Proof:\\
By [6] we have the following exact sequence:\\
$$0\rightarrow {\varprojlim}^1 K_{n+1}(BG_i)\rightarrow K_n(BG) \rightarrow \varprojlim K_n(BG_i)\rightarrow 0$$
Let us show that ${\varprojlim}^1 K_n(BG_i)=0.$ for any $n>0.$\\
We prove that the sequence $K_n(BG_i)$ is a direct summand of the sequence $K_n(BB_i).$\\
By proposition 1 of [2] We have $K_n(BG_i)=K_n^G(EG_i)$
Since we can choose $EG_i$ as a model for $EB_i,$ we obtain
$K_n(BB_i)=K_n^B(EB_i)=K_n^B(EG_i)=K_n^G(EG_i\times G/B).$\\
So, in fact, we prove that the sequence $K_n^G(EG_i)$ is a direct summand of the sequence $K_n^G(EG_i\times G/B).$\\
To simplify the notation denote $\mathcal{P}_j=\mathcal{P}_{\pi_{EG_j}}(G;EG_j\times G/B).$ 
By lemmas 6 and 7 we obtain a commutative diagram with exact arrows:
\begin{equation}
\xymatrix{
\Pcat(G;EG_j)\ar[d]_{\pi_j^*} && \Pcat(G;EG_{j+1})\ar[ll]^{(i_j\times id)^*}\ar[d]_{\pi_{j+1}^**}\\
\Pcat_j\ar[d]_{\pi_{j *}} && \Pcat_{j+1}\ar[ll]^{(i_j\times id)^*}\ar[d]_{\pi_{j *}}\\
\Pcat(G;EG_j) && \Pcat(G;EG_{j+1})\ar[ll]^{i_j^*}\\
}\label{Categories diagram}
\end{equation}

By lemma 8 the composition
$$\xymatrix{
\Pcat(G;EG_j)\ar[r]^<<<<<<{\pi_j^*} & \Pcat_j\ar[r]^<<<<<<{\pi_{j *}} & \Pcat(G;EG_j)\\
}$$
is naturally isomorphic to $id_{\Pcat(G;EG_j)}.$
In the proof of Lemma 6 it is checked that $(i_j\times id)^*(\Pcat_{j+1})\subseteq\Pcat_j$
By Lemma 5, each G-module in $\Pcat(G;EG_j\times G/B)$ has a finite resolution consisting of sheaves from $\Pcat_j.$
Then by the Quillen's theorem we get the isomorphisms $\alpha_j$ such that the following diagram of groups commutes:
\begin{equation}
\xymatrix{
K_n(\Pcat_j)\ar[d]_{\alpha_j} && K_n(\Pcat_{j+1})\ar[d]_{\alpha_{j+1}}\ar[ll]^{(i_j\times id)^*}\\
K_n^G(EG_j\times G/B) && K_n^G(EG_{j+1}\times G/B)\ar[ll]^{(i_j\times id)^*}\\
}\label{Quillen's diagram}
\end{equation}
In Corollary 3 we defined $\pi_{j*}:K_n^G(EG_j\times G/B)\rightarrow K_1^G(EG_j)$ as the composition of
$$
\xymatrix{
K_1^G(EG_j\times G/B)\ar[rr]^{\alpha_j^{-1}} && K_1(\Pcat_j)\ar[rr]^{\pi_{j*}} && K_1^G(EG_j)\\
}
$$
Commutativity of the diagrams (\ref{Categories diagram}) and (\ref{Quillen's diagram}) gives us a commutative diagram:
\begin{equation}
\xymatrix{
K_n^G(EG_j)\ar[d]_{\pi_j^*} && K_n^G(EG_{j+1})\ar[ll]^{(i_j\times id)^*}\ar[d]_{\pi_{j+1}^**}\\
K_n^G(EG_j\times G/B)\ar[d]_{\Pi_{j *}} && K_n^G(EG_{j+1}\times G/B)\ar[ll]^{(i_j\times id)^*}\ar[d]_{\pi_{j+1 *}}\\
K_n^G(EG_j) && K_n^G(EG_{j+1})\ar[ll]^{i_j^*}\\
}
\end{equation}
As we have shown, compositions of vertical arrows are identity, so $K_n^G(EG_j)$ is a direct summand of sequence $K_n^G(EG_i\times G/B)=K_n(BB_j).$
Since $\varprojlim^1(K_n(BB_j))=0$ we get $\varprojlim^1(K_n^G(EG_j))=0.$
It remains us to show that $\varprojlim^1(K_n(BB_j))=0.$ Let $T$ be a maximal torus. 
Since $B/T$ is affine space, we have that $BT_j\rightarrow BB_j$ is locally trivial with 
strats being affine spaces. Then pullback map $K_n(BB_j)\rightarrow K_n(BT_j)$ is a natural isomorphism. Since $G$ is split, $T$ is a split torus, $T=\Gm\times\ldots\times\Gm.$ Then $BT_j=\Pt^j\times\ldots\times\Pt^j.$ So, $K_n(BT_j)=K_n(pt)[t_1\ldots t_n]/(t_1^{j+1},\ldots t_n^{j+1}).$ Embedding 
pullbacks act as follows:
$$t_k \text{ mod } (t_1^{j+1},\ldots t_n^{j+1})\mapsto t_k \text{ mod } (t_1^{j},\ldots t_n^{j})$$ 
Then all morphisms in the sequence $\ldots\rightarrow K_n(BT_{j})\rightarrow K_n(BT_{j-1})\rightarrow\ldots$ are surjective. Then $\varprojlim^1(K_n(BT_i))=0,$ and consequently, $\varprojlim^1(K_n(BB_i))=0.$
 This concludes the proof.

\section{proof of main result}

{\bf{Theorem 2.}} The Borel construction induces an isomorphism
$$\xymatrix{\widehat{K_n^B(pt)}_{I_B}\ar[rr]^{\widehat{Borel^B_n}} && \widehat{K_n(BB)}_{I_B} && K_n(BB)\ar[ll]_{\cong} \\}$$

Proof:\\
We define $Borel^B_n:K_n^B(pt)\rightarrow K_n(BB)$ in the following way:
For any $j$ we construct $(Borel^B_n)_j:K_n^B(pt)\rightarrow K_n^B(EB_j)$ as a pullback of a projection $\pi_{pt}:EB_j\rightarrow pt.$
By proposition 1[2] $K_n^B(EB_j)$ are isomorphic to $K_n(BB_j).$
So we get $(Borel^B_n)_j:K_n^B(pt)\rightarrow K_n(BB_j)$
By propostion 3, we obtain $Borel^B_n=\varprojlim (Borel^B_n)_j:K_n^B(pt)\rightarrow K_n(BB).$

Let $T$ be a maximal torus of $G.$ By Corollary 1 of [2] exact functor $\Pcat(T;pt)\rightarrow \Pcat(B;B/T)$
induces an isomorphism $K_n^T(pt)\cong K_n^B(B/T).$
Since $B/T$ is affine space, we have by theorem 3 of [2] that the pullback morphism $K_n^B(pt)\rightarrow K_n^B(B/T)$ is an isomorphism.
Recall that $BT_j\rightarrow BB_j$ is locally trivial with 
strats being affine spaces. Then by theorem 3 of [2] pullback map $K_n(BB_j)\rightarrow K_n(BT_j)$ is an isomorphism. 
So we get the commutative diagram 
$$
\xymatrix{
K_n^B(pt)\ar[d]^{\cong}\ar[rr]^{Borel_n^B} && K_n(BB)\ar[d]^{\cong} \\
K_n^B(B/T)\ar[d]^{\cong}\ar[rr]^{\pi_{B/T}^*} && K_n(BB\times B/T)\ar[d]^{\cong} \\
K_n^T(pt)\ar[rr]^{Borel_n^T} && K_n(BT) \\
}
$$
So, it suffices to prove our theorem for maximal torus $T.$
Since $G$ is split,\\
$T=\Gt_m\times\ldots\times \Gt_m$ (j times).
Let us compute $K_n^T(pt)$ and $K_n^T(pt)_{I_T}.$\\
It is known that $K_n^T(pt)=K_n(pt)\otimes_{\Zt}R(T).$
$$R(T)=\Zt[\lambda_1\ldots\lambda_j,t]/(\lambda_1\cdot\ldots\cdot\lambda_j\cdot t=1).$$
$I_T=(1-\lambda_1,\ldots,1-\lambda_j,1-t).$
So, we have:
$$\widehat{K_n^T(pt)_{I_T}}=\widehat{R(T)}_{I_T}\otimes_{\Zt}K_n(pt)$$
$$\widehat{R(T)}_{I_T}=\varprojlim \Zt[\lambda_1,\ldots,\lambda_j,t]/((\Pi\lambda_i\cdot t-1), (1-\lambda_1)^k,\ldots,(1-\lambda_j)^k,(1-t)^k)=$$
$$=\varprojlim \Zt[1-\lambda_1,\ldots,1-\lambda_j,1-t]/((\Pi\lambda_i\cdot t-1), (1-\lambda_1)^k,\ldots,(1-\lambda_j)^k,(1-t)^k)=
$$
$$
=\Zt[[1-\lambda_1,\ldots,1-\lambda_j,1-t]]/(\Pi\lambda_i\cdot t-1)=
\Zt[[\mu_1,\ldots,\mu_l,1-t]]/(\Pi(1-\mu_i)\cdot t-1)
$$
Since $\frac{1}{1-\mu_i}=1+\mu_i+\mu_i^2+\mu_i^3+\ldots $  it follows that $t=\prod(1+\mu_i+\mu_i^2+\ldots)$
Therefore we have $1-t=1-(1+\mu_1 +\ldots+\mu_j+\ldots)=-(\mu_1+\ldots+\mu_j+\ldots).$ Then
$$
\widehat{R(T)}_{I_T}=\Zt[[\mu_1,\ldots,\mu_j]].
$$
Finally, we get $$\widehat{K_n^T(pt)_{I_T}}=K_n(pt)[[\mu_1,\ldots,\mu_j]]$$
Let us compute $K_n(BT).$\\
We can choose by $ET$ the space $\mathbb{A}^{\infty}\backslash \{0\}\times\ldots\times\mathbb{A}^{\infty}\backslash \{0\}$ This is contractible space with free $T-$action. Then
$ET_k=\mathbb{A}^{k+1}\backslash \{0\}\times\ldots \mathbb{A}^{k+1}\backslash \{0\}$ and
$BT_k=\Pt^{k}\times\ldots\times\Pt^{k}.$ Then $K_n(BT_k)=K_n(pt)[x_1\ldots x_n]/(x_1^k,\ldots,x_n^k).$\\
So we have $BT=\Pt^{\infty}\times\ldots\times\Pt^{\infty},$  And finally we get
$$K_n(BT)=\varprojlim K_n(BT_k)=K_n(pt)[[x_1\ldots x_n]].$$

Borel construction $K_n^T(pt)\rightarrow K_n(BT_k)$ works as follows:\\
$\lambda_i\mapsto 1-x_i$\\
$t\mapsto \frac{1}{(1-x_1)\ldots(1-x_n)}=(1+x_1+\ldots+x_1^{k-1})\ldots(1+x_1+\ldots+x_1^{k-1})$\\
Then on $\widehat{K_n^T(pt)}_{I_T}.$ Borel construction induces an isomorphism $\mu_i\mapsto x_i.$
Let us prove that $K_n(BT)$ is complete in the $I_T$-adic topology. $R(T)$-module structure on $K_n(BT)$ arises
from $R(T)$-structure on $K_0(BT)=\Zt[[x_1\ldots x_n]].$ Then $I_T\cdot K_n(BT)=(x_1,\ldots,x_n).$ So, $K_n(BT)$ is
complete.
This completes the proof of theorem 2.

\eject
\begin{flushleft}
{\bf{Theorem 3.}} There is a commutative diagram of the form:
\end{flushleft}
\noindent
\begin{equation}
\xymatrix{
     \widehat{K_n^G(pt)}_{I_G}\ar[rr]^{\widehat{Borel^G_n}}\ar[d]_{\alpha} && \widehat{K_n(BG)}_{I_G}\ar[d]_{\widehat{p^*}} && K_n(BG)\ar[ll]_{completion_G}\ar[d]_{p^*} \\
     \widehat{K_n^B(pt)}_{I_B}\ar[rr]^{\widehat{Borel^B_n}}\ar[d]_{\beta} && \widehat{K_n(BB)}_{I_B}\ar[d]_{\widehat{p_*}} && K_n(BB)\ar[ll]_{completion_B}\ar[d]_{p_*} \\
     \widehat{K_n^G(pt)}_{I_G}\ar[rr]^{\widehat{Borel^G_n}} && \widehat{K_n(BG)}_{I_G} && K_n(BG) \ar[ll]_{completion_G} \\
}
\label{bebe main reduction diagram}
\end{equation}

with $\beta\circ\alpha=id,$ $\widehat{p_*}\circ\widehat{p^*}=id,$ and  $p_*\circ p^*=id.$ \\
Proof:\\
Since $EG_i\rightarrow BG_i$ is a $G$-torsor, $K_n(BG_i)=K_n^G(EG_i).$ (by Proposition 1 of[2])
$EG$ can be chosen as a model for the contractible space $EB$
Proposition 1 of [2] allows us express all these objects in terms of $G$-equivariant K-theory:
$K_n^B(pt)\cong K_n^G(G/B)$ $K_n^B(EG_j)=K_n^G(EG_j\times G/B)$

So, first we construct :
\begin{equation}
\xymatrix{
     K_n^G(pt)\ar[rr]^{\pi_{pt}^*}\ar[d]_{\pi_{pt}^*} && K_n^G(EG_i)\ar[d]_{\pi_{EG_i}^*} \\
     K_n^G(G/B)\ar[rr]^{\pi_{G/B}^*}\ar[d]_{\pi_{pt *}} && K_n^G(EG_i\times G/B)\ar[d]_{\pi_{EG_i *}} \\
     K_n^G(pt)\ar[rr]^{\pi_{pt}^*} && K_n^G(EG_i) \\
}
\label{j level reduction diagram}
\end{equation}
Proposition 1 proves that this diagram commutes and 
$\pi_{pt*}\circ \pi_{pt}^*=id$ and $\pi_{EG_i *}\circ \pi_{EG_i}^*=id.$
Recall that $K_n^G(EG_j)=K_n(BG_j)$ , $K_n^G(EG_j\times G/B)=K_n(BB_j),$ and $K_n^G(G/B)=K_n^B(pt).$\\
Therefore we can rewrite the above diagram as follows
\begin{equation}
\xymatrix{
     K_n^G(pt)\ar[rr]^{\pi_{pt}^*}\ar[d]_{\pi_{pt}^*} && K_n(BG_i)\ar[d]_{\pi_{EG_i}^*} \\
     K_n^B(pt)\ar[rr]^{\pi_{G/B}^*}\ar[d]_{\pi_{pt *}} && K_n(BB_i)\ar[d]_{\pi_{EG_i *}} \\
     K_n^G(pt)\ar[rr]^{\pi_{pt}^*} && K_n(BG_i) \\
}
\label{j level nonequivariant reduction diagram}
\end{equation}
Take the projective limit of this diagram. Recall that $K_n(BB_i)=K_n(BB)$ and by proposition 3 we have $\varprojlim K_n(BG_i)=K_n(BG).$
So we get commutative diagram of $K_0^G(pt)$-modules
\begin{equation}
\xymatrix{
     K_n^G(pt)\ar[rr]^{Borel_n^G}\ar[d]_{\pi_{pt}^*} && K_n(BG)\ar[d]_{\varprojlim\pi_{EG_i}^*} \\
     K_n^B(pt)\ar[rr]^{Borel_n^B}\ar[d]_{\pi_{pt *}} && K_n(BB)\ar[d]_{\varprojlim\pi_{EG_i *}} \\
     K_n^G(pt)\ar[rr]^{Borel_n^G} && K_n(BG) \\
}
\label{final reduction diagram}
\end{equation}
Here we still have $\pi_{pt *}\circ\pi_{pt}^*=id$ and $\varprojlim\pi_{EG_i *}\circ\varprojlim\pi_{EG_i}^*=id.$
Let us denote $p_*=\varprojlim\pi_{EG_i *}$ and $p^*=\varprojlim\pi_{EG_i}^*.$
Recall that $R(G)$-structures on $K_n(BB)$ and $K_n^B(pt)$ are induced by
$R(G)$-structure on $R(B).$ Then proposition 2 implies that $I_G$-adic  completions of $K_n(BB)$ and $K_n^B(pt)$ coincides with
$I_B$-adic completions. So, by taking $I_G$-adic completion of(\ref{final reduction diagram}) we obtain commutative diagram
\begin{equation}
\xymatrix{
     \widehat{K_n^G(pt)}_{I_G}\ar[rr]^{\widehat{Borel_n^G}}\ar[d]_{\widehat{\pi_{pt}^*}} && \widehat{K_n(BG)}_{I_G}\ar[d]_{\widehat{p^*}} \\
     \widehat{K_n^B(pt)}_{I_B}\ar[rr]^{\widehat{Borel_n^B}}\ar[d]_{\widehat{\pi_{pt *}}} && \widehat{K_n(BB)}_{I_B}\ar[d]_{\widehat{p_*}} \\
     \widehat{K_n^G(pt)}_{I_G}\ar[rr]^{\widehat{Borel_n^G}} && \widehat{K_n(BG)}_{I_G} \\
}
\label{completed reduction diagram}
\end{equation}
with $\widehat{\pi_{*}}\circ\widehat{\pi^*}=id$ and $\widehat{p_{*}}\circ \widehat{p^*}=id.$
Consider the commutative diagram:
\begin{equation}
\xymatrix{
     \widehat{K_n(BG)}_{I_G}\ar[d]_{\widehat{p^*}} && {K_n(BG)}\ar[d]_{p^*}\ar[ll]_{completion_G}\\
     \widehat{K_n(BB)}_{I_B}\ar[d]_{\widehat{p_*}} && {K_n(BB)}\ar[d]_{p_*}\ar[ll]_{completion_B}\\
     \widehat{K_n(BG)}_{I_G} && {K_n(BG)}\ar[ll]_{completion_G}\\
}
\label{obvious}
\end{equation}
Set $\alpha=\widehat{\pi_{pt}^*},$ $\beta=\widehat{\pi_{pt *}}$ and recall that $K_n^G(G/B)=K_n^B(pt).$ Then by gluing together $\ref{obvious}$ and $\ref{completed reduction diagram}$ we obtain the diagram $\ref{bebe main reduction diagram}:$
$$
\xymatrix{
     \widehat{K_n^G(pt)}_{I_G}\ar[rr]^{\widehat{Borel^G_n}}\ar[d]_{\alpha} && \widehat{K_n(BG)}_{I_G}\ar[d]_{\widehat{p^*}} && K_n(BG)\ar[ll]_{completion_G}\ar[d]_{p^*} \\
     \widehat{K_n^B(pt)}_{I_B}\ar[rr]^{\widehat{Borel^B_n}}\ar[d]_{\beta} && \widehat{K_n(BB)}_{I_B}\ar[d]_{\widehat{p_*}} && K_n(BB)\ar[ll]_{completion_B}\ar[d]_{p_*} \\
     \widehat{K_n^G(pt)}_{I_G}\ar[rr]^{\widehat{Borel^G_n}} && \widehat{K_n(BG)}_{I_G} && K_n(BG) \ar[ll]_{completion_G} \\
}
$$

with $\beta\circ\alpha=id,$ $\widehat{p_*}\circ\widehat{p^*}=id,$ and  $p_*\circ p^*=id.$ This concludes the proof.\\

{\bf{Theorem 1}} In the following diagram both maps are $K_0^G(pt)$-module isomorphisms:
$$\xymatrix{\widehat{K_n^G(pt)_{I_G}}\ar[rr]^{\widehat{Borel^G_n}} && \widehat{K_n(BG)}_{I_G} && K_n(BG)\ar[ll]_{completion_G} \\}$$
Proof:\\
Theorem 3 states that $\widehat{Borel^G_n}$ and $completion_G$
are retracts of $\widehat{Borel^B_n}$ and $completion_B$ which are isomorphisms by theorem 2. Then $\widehat{Borel^G_n}$ and $completion_G$ are also isomorphisms.

 \vfil

\centerline{References:}
\noindent 1. R. Hartshorne. {\it Algebraic Geometry}. New York, Springer-Verlag(1977).\\
2. A.S. Merkurjev. {\it Equivariant K-theory}. Handbook of K-theory. Vol 1, 2:925-954, Springer(2005).\\
3. R.W. Thomason. {\it Algebraic K-theory of group scheme actions}, Algebraic topology and algebraic K-theory 
(Princeton, N.J., 1983), Ann. of Math. Stud., vol.113, Princeton Univ.Press, Princeton, N.J, pp.539-563,(1987).\\
4. D. Mumford. {\it Abelian Varieties}. Oxford University Press(1988).\\
5. D. Mumford, J. Fogarty, and F. Kirwan. {\it Geometric invariant theory}. New York, Springer-Verlag(1994).\\
6. F. Morel, V. Voevodsky. {\it $\mathbb{A}^1$-homotopy theory of schemes}. Publ. Math., Inst. Hautes \'Etud. Sci. 90: 45-143, (1999).\\
7. J.С. Jantzen. {\it Representations of algebraic groups}. Academic press, Orlando-London(1987).\\
8. B. Totaro. {\it The Chow ring of a classifying spaces}. "`Algebraic K-theory"'(Seattle, WA, 1997):249-281, Proc. Sympos. Pure Math. 67, Amer Math. Soc.(1999).\\
9. D. Quillen. {\it Higher algebraic K-theory I}. Lecture Notes in Math., Vol.341:85-147, (1973).\\
10. V. Voevodsky. {\it $\mathbb{A}^1$-homotopy theory}. Documenta Mathematica, Extra Volume ICM 1998(I):579-604, (1998).

\end{document}